\documentclass[11pt,fleqn,twoside]{article}
\usepackage{emlines,latexsym,amsfonts}
\usepackage[dvips]{epsfig}
\makeatletter
\newcommand{\prava}{\footnotesize\it
\begin{flushright}
\begin{minipage}{18cm}%{6cm}%9.6
Copyright \copyright 1998 by M.B. Abd-el-Malek
\end{minipage}
\end{flushright}}

\newcommand{\name}[1]{\begin{flushleft}
                       \LARGE \bf #1
                       \end{flushleft}\vspace{-3mm}}

\newcommand{\Author}[1]{\begin{flushleft}
                       \it #1 \end{flushleft}}

\newcommand{\Adress}[1]{\begin{flushleft}
                       \it #1 \end{flushleft}}

\newcommand{\Date}[1]{\begin{flushleft}
                      \small  \it #1 \end{flushleft}}

\newcommand{\ehkol}{Author \ name}
\newcommand{\ohkol}{Article \ name}
\renewcommand{\@evenhead}{
\hspace*{-3pt}\raisebox{-15pt}[\headheight][0pt]{\vbox{\hbox to \textwidth
{\thepage \hfil \ehkol}\vskip4pt \hrule}}}
\renewcommand{\@oddhead}{
\hspace*{-3pt}\raisebox{-15pt}[\headheight][0pt]{\vbox{\hbox to \textwidth
{\ohkol \hfil \thepage}\vskip4pt\hrule}}}
\renewcommand{\@evenfoot}{}
\renewcommand{\@oddfoot}{}

\newcommand{\be}{\begin{equation}}
\newcommand{\ee}{\end{equation}}
\newcommand{\ba}{\hspace*{-5pt}\begin{array}}
\newcommand{\ea}{\end{array}}
\newcommand{\p}{\partial}
\newcommand{\ds}{\displaystyle}
\makeatother

\begin{document}
\setcounter{page}{314}
\thispagestyle{empty}

\renewcommand{\ehkol}{M.B. Abd-el-Malek}
\renewcommand{\ohkol}{Application of the Group-Theoretical Method}

\begin{flushleft}
\footnotesize \sf
Journal of Nonlinear Mathematical Physics \qquad 1998, V.5, N~3,
\pageref{malek-fp}--\pageref{malek-lp}.\hfill {\sc Article}
\end{flushleft}

\vspace{-5mm}

\renewcommand{\footnoterule}{}
{\renewcommand{\thefootnote}{}  
\footnote{\prava}}

\name{Application of the Group-Theoretical Method \\
to Physical Problems}\label{malek-fp} 

\Author{Mina B. ABD-EL-MALEK}

\Adress{Department of Engineering Mathematics and Physics, Faculty of
Engineering,\\
Alexandria University, Alexandria 21544, Egypt}

\Date{Received 9 March, 1998}

\begin{abstract}
\noindent
The concept of the theory of continuous groups of transformations has
attracted the attention of applied mathematicians and engineers to
solve many physical problems in the engineering sciences. Three
applications are presented in this paper. The f\/irst one is the
problem of time-dependent vertical temperature distribution in a
stagnant lake. Two cases have been considered for the forms of the
water parameters, namely water density and thermal conductivity. The
second application is the unsteady free-convective boundary-layer
f\/low on a non-isothermal vertical f\/lat plate. The third
application is the
study of the dispersion of gaseous pollutants in the presence of a
temperature inversion. The results are found in closed form and the
ef\/fect of parameters are discussed.
\end{abstract}

\section{Introduction}

According to 
Seshadri and Na [1],
dif\/ferent methods for carrying out similarity analysis of partial
dif\/ferential equations are classif\/ied 
into direct methods, where the concept of group
invariance is not explicitly invoked, and 
group-theoretical methods that are based upon the invocation of invariance
under groups of transformations of the 
partial dif\/ferential equation and the auxiliary conditions.

The foundation of the group-theoretical method is contained in the general
theory of continuous transformation 
groups that were introduced and treated extensively by Lie [2] in 1875.
Group-theoretical methods provide a powerful 
tool because they are not based on linear operators, superposition, or any
other aspects of linear solution 
techniques. Therefore, these methods are applicable to nonlinear
dif\/ferential models. 

Throughout the history of similarity analysis, a variety of problems in
science and engineering have been solved. 
Among these we f\/ind a general treatment of steady two-dimensional
incompressible laminar f\/low of f\/luid into an 
inf\/inite region of the same f\/luid by Abbott
and Kline [3] in 1960, impact of
thin long rods by Taulbee {\it et al} [4] in 
1971, wave propagation in viscoelastic, viscoplastic and electrical
transmission lines by Ames and Suliciu [5] in 
1982, time-dependent free surface f\/lows under gravity by Sachdev and Philip
[6] in 1986, problem of ocean 
acoustics by Richards [7] in 1987, unsteady free-convective boundary-layer
f\/low on a non-isothermal vertical f\/lat 
plate by Abd-el-Malek {\it et al} [8] in 1990, steady free-convective
boundary-layer f\/low on a non-isothermal vertical 
circular cylinder by Abd-el-Malek and Badran [9] in 1991, dispersion of
gaseous pollutants in the presence of a 
temperature inversion by Badran and Abd-el-Malek [10] in 1993, and
nonlinear temperature variation across the lake 
depth neglecting the ef\/fect of external heat sources by Abd-el-Malek [11] in
1997. Recently Boutros {\it et al} [12, 13] 
considered two problems namely, potential equation in triangular regions
with temperature distribution along the 
boundaries in the form of polynomials, and the second is the
time-dependent vertical temperature distribution 
in stagnant lake taking into consideration an external heat source. Many
physical application are illustrated in Sedov [14], and Rogers and Ames~[15].

In the present paper we will present three physical applications of the
group-theoretical method namely, time-dependent vertical temperature
distribution in a stagnant lake taking into consideration the ef\/fect of
external heat source. The problem has been solved for two possible forms of
the water parameters (water density and thermal 
conductivity). Second application is the unsteady free-convective
boundary-layer f\/low on a non-isothermal 
vertical f\/lat plate. As a third application we study the dispersion of
gaseous pollutants in the presence of a 
temperature inversion. The obtained results are found in closed form and the
ef\/fect of dif\/ferent parameters are discussed.

\renewcommand{\thesection}{\arabic{section}}
\setcounter{equation}{0}
\renewcommand{\theequation}{\arabic{section}.\arabic{equation}}

\section{Application (I): Time-dependent vertical temperature\\
distribution in a stagnant lake}

\subsection{Mathematical formulation}

Consider the one dimensional heat transfer equation for heat f\/lux in the
vertical direction, neglecting the convective motion of the f\/luid and
assuming the absolute value of the specif\/ic heat of the water as sensibly 
constant within the range of temperatures considered; take it to be
unity.
The
vertical transport of heat in a deep lake is 
modelled by
\be \label{malek:I.1}
\rho(T)T_t =[\kappa(T)T_z]_z +r(z,t), \qquad z>0, \quad t>0,
\ee
where ``$T$'' is the temperature; $r(z,t)$ 
is rate at which solar radiation is absorbed by the water; ``$t$'' 
is the time; ``$z$'' is the distance measured downward from the lake
surface; ``$\rho$'' is the density; and ``$\kappa$'' 
is the thermal conductivity. 

\medskip

\noindent
{\bf Boundary and initial conditions:}

\setcounter{equation}{0}
\renewcommand{\theequation}{\arabic{section}.2.\arabic{equation}}

\noindent
During early Spring, most of the lakes exhibit a nearly homothermal
temperature distribution with a low degree 
of temperature extending all the way to the bottom (see Sundaram and Rehm
[16]). In all of the calculations 
presented in this paper, the initial condition will be taken as that
corresponding to the end of spring homothermy, i.e.,
\be \label{malek:I.2.1}
T(z,0) = T_0, 
\ee
where $T_0$ is the temperature of the lake at maximum Spring homothermy.

{\it The boundary conditions} are considered as follows:
\be\label{malek:I.2.2}
\mbox{(i)}\  \quad T_z(0,t)=0, \quad t>0,
\ee
\be \label{malek:I.2.3}
\mbox{(ii)} \quad T(h,t)=T_0+\gamma t^{1/m}, \quad t>0, \ 0\leq \gamma \ll 1,
\ m>0.
\ee

In our analysis we will consider two cases for the form of the density
$\rho(T)$ and the thermal conductivity $\kappa(T)$, 
namely:
\[
\mbox{Case (1):} \quad 
\rho=\alpha q(z)(T-T_0)^m, \quad \kappa=\beta g(z),
\]
\[
\mbox{Case (2):} \quad 
\rho=\alpha q(z)(T-T_0)^s, \quad \kappa=\beta(T-T_0)^n,
\]
where $\alpha$ and $\beta$ are constants, $m$, $s$ and $n$ are positive
constants and $q(z)$ and $g(z)$ are arbitrary functions to be 
determined later on.

\setcounter{equation}{2}
\renewcommand{\theequation}{\arabic{section}.\arabic{equation}}

Write
\be \label{malek:I.3}
T(z,t)=w(z,t)+T_0,
\ee
by which dif\/ferential equation (\ref{malek:I.1}) takes the form:
\be \label{malek:I.4}
\rho(w) w_t =[k(w)w_z]_z+r(z,t), \quad z>0, \ t>0,
\ee
and the initial and boundary conditions take the form:

\setcounter{equation}{0}
\renewcommand{\theequation}{\arabic{section}.5.\arabic{equation}}

(1) {\it initial condition}:
\be \label{malek:I.5.1}
w(z,0)=0,
\ee

(2) {\it boundary conditions}:
\be \label{malek:I.5.2}
\mbox{(i)} \ \quad w_z(0,t)=0, \quad t>0,
\ee
\be \label{malek:I.5.3}
\mbox{(ii)} \quad w(h,t)=\gamma t^{1/m}, \quad t >0, \
0\leq \gamma \ll 1, \ m>0.
\ee

\subsection{Case (1): $\rho=\alpha q(z)w^m$, $\kappa=\beta g(z)$}

\setcounter{equation}{5}
\renewcommand{\theequation}{\arabic{section}.\arabic{equation}}

In this case dif\/ferential equation (\ref{malek:I.4}) takes the form:
\be \label{malek:I.6}
\beta g(z) w_{zz} +\beta g_z w_z -\alpha q(z) w^m w_t=
- r(z,t).
\ee

\subsubsection{Solution of the problem}

Our method of solution depends on the application of a one-parameter group
transformation to the partial dif\/ferential equation (\ref{malek:I.6}).

\smallskip

\noindent
{\bfseries \itshape 2.2.1 (a) The group systematic formulation.}
The procedure is initiated with the group $G$; a class of transformation of
one-parameter ``$a$" of the form 
\be \label{malek:I.7}
G: \quad \overline{Q} =C^Q(a)Q+P^Q(a),
\ee
where $Q$ stands for $t$, $z$, $r$, $w$, $\kappa$, $\rho$ and 
the $C$'s and $P$'s are real-valued and at least dif\/ferentiable in ``$a$''.

\smallskip

\noindent
{\bfseries \itshape 2.2.1 (b) The invariance analysis.}
To transform the dif\/ferential equation, transformations of the derivatives
of $w$, $\kappa$ and $\rho$ are obtained from $G$ via 
chain-rule operations:
\be \label{malek:I.8}
\overline{Q}_{\overline{i}} =\left(\frac{C^Q}{C^i}\right) Q_i,
\qquad \overline{Q}_{\overline{i}\, \overline{j}}=
\left(\frac{C^Q}{C^i C^j}\right) Q_{ij}; 
\qquad i=z,t; \ j=z,t.
\ee
Equation (\ref{malek:I.6}) is said to be invariantly transformed, for some
function $H(a)$, whenever 
\be \label{malek:I.9}
\beta \overline{g} \, \overline{w}_{\overline{z}\, \overline{z}}+
\beta \overline{g}_{\overline{z}}\, \overline{w}_{\overline{z}}-
\alpha \overline{q}(\overline{w})^m \,\overline{w}_{\overline{t}} +
\overline{r} =H(a) \left[ \beta g w_{zz} +\beta g_z w_z -\alpha q w^m w_t+r
\right] .
\ee
 Substitution from (\ref{malek:I.7}) into (\ref{malek:I.9}) yields
\be \label{malek:I.10}
\ba{l}
\ds \beta \left[ \frac{C^g C^w}{(C^z)^2}\right] g w_{zz} +
\beta\left[\frac{C^g C^w}{(C^z)^2}\right] g_z w_z -
\alpha w^m \left[\frac{C^q(C^w)^{m+1}}{C^t}\right] q w_t +
[C^r]r +R(a)\\[5mm]
\ds \qquad = H(a)\left[ \beta g w_{zz} +\beta g_z w_z -\alpha q w^m w_t+r
\right],
\ea
\ee
where
\[
\ba{l}
\ds R(a) = \left[ \beta \frac{P^g C^w}{(C^z)^2}\right] w_{zz}
-\left[ \alpha P^q(C^w w+P^w)^m \frac{C^w}{C^t}\right]w_t\\[5mm]
\ds \qquad - \alpha(C^q q) \frac{C^w}{C^t} w_t
\sum\limits_{k=1}^m \left( \ba{c} m \\ k \ea \right)
(C^w w)^{m-k}(P^w)^k +P^r.
\ea
\]
The invariance of (\ref{malek:I.9}) implies $R(a)\equiv 0$. This is satisf\/ied by
putting 
\be \label{malek:I.11}
P^r=P^q=P^w=P^g=0,
\ee
and
\be \label{malek:I.12}
\left[\frac{C^q (C^w)^{m+1}}{C^t}\right] =
\left[\frac{C^g C^w}{(C^z)^2}\right] = C^r = H(a).
\ee

Moreover, the boundary and initial conditions (\ref{malek:I.5.2}), (\ref{malek:I.5.3})
and (\ref{malek:I.5.1}) are also invariant in form, implying that 
\be \label{malek:I.13}
P^z=0 , \quad P^T=0, \quad C^z=1 \quad \mbox{and} \quad  (C^w)^m = C^t. 
\ee
Combining equations (\ref{malek:I.12}) and invoking the result (\ref{malek:I.13}), we
get 
\be \label{malek:I.14}
C^r=C^qC^w =C^g C^w \qquad \mbox{which yields} \quad C^q=C^g.
\ee

Finally, we get the one-parameter group G which transforms invariantly the
dif\/ferential equation (\ref{malek:I.1}), as well as the 
boundary and initial conditions (2.2). The group $G$ is of the
form 
\be \label{malek:I.15}
G: \quad \overline{z}=z, \  \
\overline{t}=(C^w)^m t , \ \ \overline{q}=
C^q q, \ \ \overline{r}=C^q C^w r, \ \ \overline{w}=C^w w, \ \
\overline{g}= C^q g.
\ee

\noindent
{\bfseries \itshape 2.2.1 (c) The complete set of absolute invariants.}
Our aim is to make use of group methods to represent the
problem in the form of an ordinary dif\/ferential
equation. Then we have to proceed in our analysis to
obtain a complete set of absolute invariants.

 If $\eta\equiv \eta(z,t)$ is the absolute invariant of the
 independent variables, then
\be \label{malek:I.16}
g_j(z,t;w,r,\kappa,\rho) = F_j[\eta(z,t)];  \qquad j=1,2,3,4
\ee
are the four absolute invariants corresponding to $w$, $r$, $\kappa$
and $\rho$. The application of a basic theorem in group
theory, see [17], states that: {\it a function $g(z,t;w,r,\kappa,\rho)$
is an absolute invariant of a one-parameter group if it
satisfies the following first-order linear differential equation}
\be\label{malek:I.17}
\sum_{i=1}^6 (\alpha_i Q_i +\beta_i)\frac{\p g}{\p Q_i}=0; \qquad
Q_i=z,t,w,r,\kappa,\rho,
\ee
where
\be \label{malek:I.18}
\alpha_i= \frac{\p C^{Q_i}}{\p a}(a^0), \qquad
\beta_i=\frac{\p P^{Q_i}}{\p a}(a^0); \qquad i=1,2,\ldots,6,
\ee
and $a^0$ denotes the value of ``$a$''
which yields the identity element of the group.

Owing to equation (\ref{malek:I.17}), $\eta(z,t)$ is an absolute
invariant if it satisf\/ies
\be \label{malek:I.19}
(\alpha_1 z + \beta_1) \eta_z + (\alpha_2 t + \beta_2) \eta_t = 0 ,
\ee
  Group (\ref{malek:I.15}) gives:
  \be \label{malek:I.20}
  \alpha_1 = \beta_1 = \beta_2 = 0,
  \ee
and hence from (\ref{malek:I.17}) and (\ref{malek:I.20}) we get
\be \label{malek:I.21}
\eta_t=0,
\ee
which gives
\be \label{malek:I.22}
\eta(z,t)=F(z).
\ee

Without loss of generality we can use the identity function:
\be \label{malek:I.23}
\eta(z,t)=z.
\ee
By a similar analysis the absolute invariants of the dependent
variables $w$, $r$, $\kappa$ and $\rho$ are
\be \label{malek:I.24}
\ba{l}
w(z,t)=\Gamma(t)F(\eta), \qquad
r(z,t)=A(t)\Theta(\eta), \\[2mm]
q(z)=B(t)\Phi(\eta), \qquad
g(z)=Y(t)\Psi(\eta).
\ea
\ee

 From which we conclude that
\be \label{malek:I.25}
q(z) = \Phi(\eta),
\ee
\be \label{malek:I.26}
g(z) = \Psi(\eta).
\ee
\be \label{malek:I.27}
\mbox{At} \ t = 0: \quad \Gamma(0) = 0.
\ee

\subsubsection{The reduction to ordinary dif\/ferential equation}

Substituting from (\ref{malek:I.24}) into (\ref{malek:I.6}) and dividing by $\Gamma$,
we get
\be \label{malek:I.28}
\beta \Psi F_{\eta\eta}+\beta \Psi_\eta F_\eta-\alpha \Phi \Gamma^{m-1}
F^{m+1}\Gamma_t=-\frac{A(t)\Theta(\eta)}{\Gamma}.
\ee

For (\ref{malek:I.28}) to be reduced to an expression in
the single independent invariant $\eta$, it is necessary that the
coef\/f\/icients should be constants or functions of $\eta$ alone. Thus
\be \label{malek:I.29}
\Gamma^{m-1}\Gamma_t=C_1, \qquad \frac{A(t)}{\Gamma}=C_2.
\ee
Take $C_1 = 1$:
\be \label{malek:I.30}
\Gamma(t)=(mt)^{1/m}, \quad m\not=0,
\ee
(\ref{malek:I.29}) and (\ref{malek:I.30}) yield:
\be \label{malek:I.31}
A(t)=C_2(mt)^{1/m}, \quad m\not= 0.
\ee
Hence (\ref{malek:I.28}) may be rewritten as
\be \label{malek:I.32}
\beta \Psi F_{\eta\eta}+\beta \Psi_\eta F_\eta -\alpha \Phi F^{m+1}=
-C_2\Theta(\eta).
\ee
Following Girgis and Smith [18], we assume the heat source
distribution in the form:
\be \label{malek:I.33}
\Theta(\eta)={\mathrm e}^{-\xi \eta},
\ee
where $\xi$ is the absorption coef\/f\/icient for water, which has the
value $0.048$.

Take:
\be \label{malek:I.34}
\Phi(\eta)=\frac{\Psi(\eta)}{F(\eta)}, \qquad
F(\eta)\not= 0, \quad 0\leq \eta \leq h.
\ee
According to (\ref{malek:I.33}) and (\ref{malek:I.34}), equation (\ref{malek:I.32})
takes the form
\be \label{malek:I.35}
\beta \Psi F_{\eta\eta}+\beta \Psi_\eta F_\eta -\alpha \Psi F^m=
-C_2{\mathrm e}^{-\xi \eta}.
\ee
Write:
\be \label{malek:I.36}
\Psi ={\mathrm e}^{-\mu\eta},
\ee
where $\mu$ is a constant, then (\ref{malek:I.35}) becomes
\be \label{malek:I.37}
F_{\eta\eta}-\mu F_\eta-\frac{\alpha}{\beta}F^m =-\frac{C_2}{\beta}
{\mathrm e}^{-(\xi-\mu)\eta},
\ee
with the boundary conditions:
\be \label{malek:I.38}
\mbox{(i)} \ \quad   F(0) = 0,
\ee
\be \label{malek:I.39}
\mbox{(ii)}  \quad F(h)  =\gamma/m^{(1/m)}.
\ee

\subsubsection{Analytical solution for dif\/ferent forms of
the parameters}

For $m = 1$, dif\/ferential equation (\ref{malek:I.37}) becomes
\be \label{malek:I.40}
F_{\eta\eta}-\mu F_\eta-\sigma^2 F =-\frac{C_2}{\beta}
{\mathrm e}^{-(\xi-\mu)\eta}; \qquad \sigma^2=\frac{\alpha}{\beta}
\ee
and the boundary conditions become
\be \label{malek:I.41}
\mbox{(i)} \ \quad   F(0) = 0,
\ee
\be \label{malek:I.42}
\mbox{(ii)}  \quad F(h)  =\gamma,
\ee
which has the exact solution
\be \label{malek:I.43}
F(\eta)=a_1{\mathrm e}^{r_1 \eta} +a_2{\mathrm e}^{r_2 \eta}
+ a_3{\mathrm e}^{-(\xi -\mu)\eta},
\ee
where 
\[
r_{1,2}= \frac{\mu\pm\sqrt{\mu^2+4\sigma^2}}{2},
\]
$r_1$
for $(+)$  sign and $r_2$ for $(-)$ sign,
and  
\[
 a_3=-\frac{C_2}{\beta(\mu-\xi)^2-\mu\beta(\mu-\xi)-\alpha}.
\]

For f\/inite temperature, $a_1 = 0$, and applying condition
(\ref{malek:I.41}), we get the solution in the form: 
\be \label{malek:I.44}
F(\eta) =a_3 \left[ {\mathrm e}^{-(\xi-\mu)\eta}+\frac{\xi-\mu}{r_2}
{\mathrm e}^{r_2\eta}\right].
\ee
Hence the temperature distribution across the lake, corresponding to
case~(1) is: 
\be \label{malek:I.45}
T(z,t) =T_0-\frac{C_2 t}{\beta(\mu -\xi)^2 -\beta
\mu(\mu-\xi)-\alpha} 
\left[ {\mathrm e}^{-(\xi-\mu)z}+\frac{\xi-\mu}{r_2}
{\mathrm e}^{r_2 z }\right].
\ee

For $0 < t < 150$ (in days), following Girgis and Smith [18], we use
the following values of the parameters: 
$C_2=2496$, $\beta = 12355$, $\mu  = 1.439239\times 10^{-4}$,
$\xi= 0.048$, $h = 400$~meter, $T_0 = 4\,{}^\circ$C. The obtained
results are plotted in Fig.1 and Fig.2.

\subsection{Case (2): $\rho=\alpha q(z) w^s$, $\kappa=\beta w^n$}

Dif\/ferential equation (\ref{malek:I.4}) takes the form
\be \label{malek:I.46}
w^n w_{zz} +n w^{n-1} (w_z)^2-\sigma^2 q(z)w^s w_t=-r(z,t).
\ee

Following the same analysis as in case (1), we get the following
group $G$: 
\be \label{malek:I.47}
G: \quad \overline{z}=z, \ \ \overline{t}=(C^w)^m t, \ \
\overline{q}=(C^w)^{n+m-s} q, \ \
\overline{r}=(C^w)^{n+1} r, \ \ 
\overline{w}=C^w w.
\ee
%%before fig12.eps%%%%%
\epsfig{file=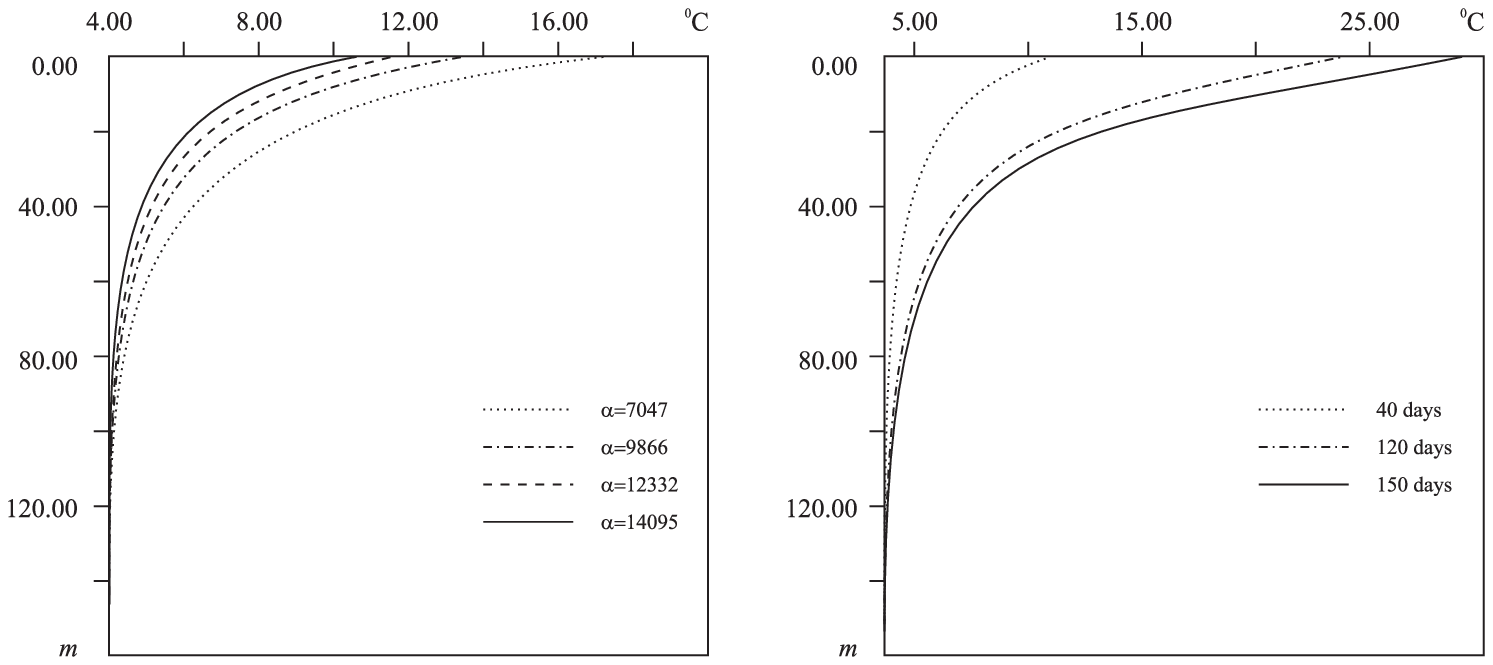,height=7cm,width=14.3cm}%14.3

%\begin{figure}[t]
%\footnotesize
%\begin{picture}(3,3)
%\put(0,0){\special{em: graph fig1-2.pcx}}%malek1.pcx
%\end{picture}

%\vspace{2mm}%68
\begin{center}
\quad\qquad
\begin{minipage}[t]{6cm}%7.1cm%6
Figure~1. Distribution of temperature $T$ (time = 40 days) againts
the lake depth ``$z$'' in meters, corresponding to case: $\rho=
\alpha q(z)w$, $\kappa=\beta g(z)$ for dif\/ferent values of
parameter ``$\alpha$''.
\end{minipage}
\hfill
\noindent
\begin{minipage}[t]{6.5cm}%
Figure~2. Distribution of temperature for dif\/ferent times to constant
``$\alpha$'' ($\alpha=14095$), for $\rho= \alpha q(z)w$, $\kappa=\beta g(z)$.
\end{minipage}
\end{center}
\vspace{3mm}%-3

%\end{figure}

The absolute invariants are:
\be \label{malek:I.48}
\eta=z, \quad w(z,t)= \Gamma(t) F(\eta), \quad
r(z,t)=A(t)\Theta(\eta),\quad q(z)=B(t)\Phi(\eta).
\ee

Again, it is clear that $B(t) = 1$; from which we get
\be \label{malek:I.49}
q(z) = \Phi(\eta),
\ee
leading to 
\be \label{malek:I.50}
\overline{q}=q,
\ee
which is satisf\/ied if and only if $(C^w)^{n+m-s}=1$.
That is $C^w = 1$ or $n = s - m$. $C^w$ can
not be unity. Hence we obtain the only possible case:
\be \label{malek:I.51}
n = s - m. 
\ee
Take
\be \label{malek:I.52}
\Phi(\eta)= \frac{1}{F(\eta)}, \qquad F(\eta)\not= 0, \quad 0\leq
\eta \leq h,
\ee
\be \label{malek:I.53}
\mbox{At} \quad  t=0: \quad \Gamma(0) = 0.
\ee
%%%put here%%%%%%%%%
%%Before figure34.eps%%%%%%%%%%%
\epsfig{file=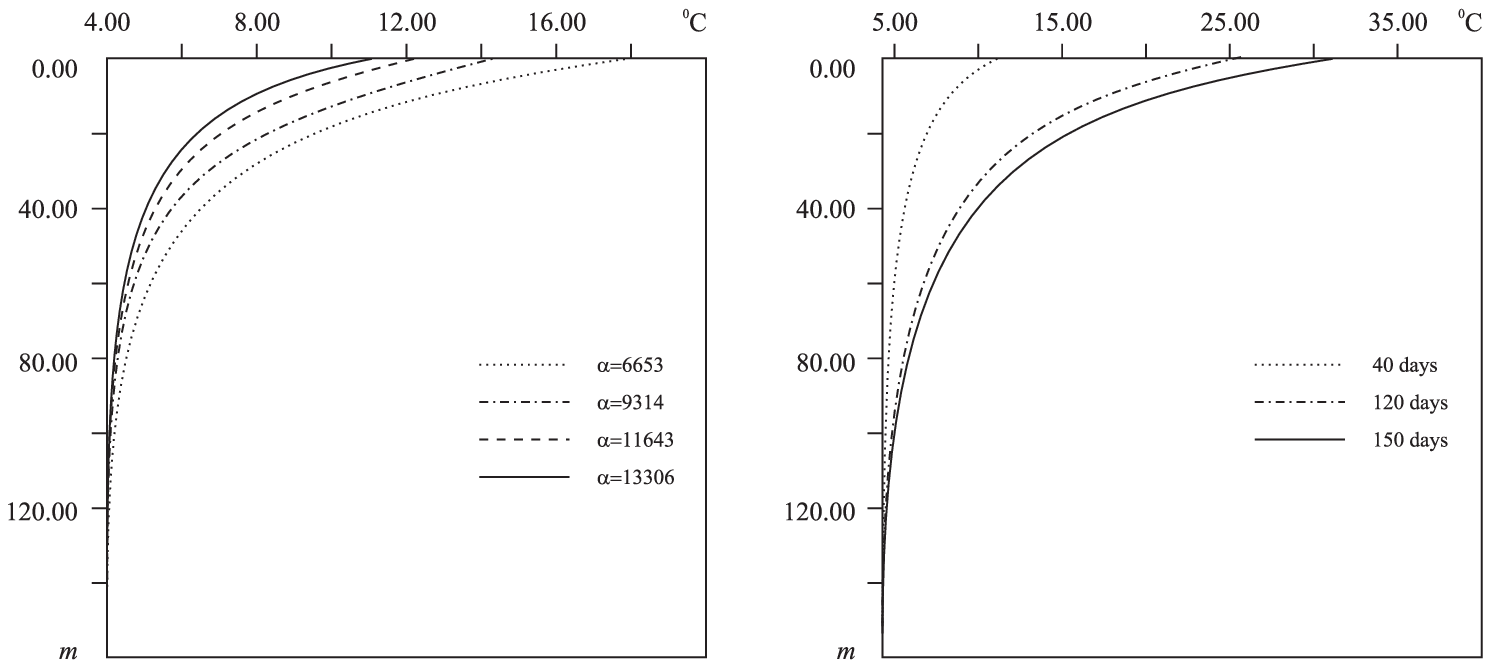,height=7cm,width=14.3cm}
%\begin{figure}[t]
%\footnotesize
%\begin{picture}(3,3)
%\put(0,0){\special{em: graph fig3-4.pcx}}%malek2.pcx
%\end{picture}

%\vspace{3mm}

\begin{center}
\quad\qquad
\begin{minipage}[t]{6cm}
Figure~3. Distribution of temperature $T$ (time = 40 days) againts
the lake depth ``$z$'' in meters, corresponding to case: $\rho=
\alpha q(z)w$, $\kappa=\beta$ for dif\/ferent values of
parameter ``$\alpha$''.
\end{minipage}
\hfill
\noindent
\begin{minipage}[t]{6.5cm}
Figure~4. Distribution of temperature for dif\/ferent times to constant
``$\alpha$'' ($\alpha=13306$), for $\rho= \alpha q(z)w$, $\kappa=\beta g(z)$.
\end{minipage}
\end{center}

\vspace{3mm}

%\end{figure}

Following the same analysis as in section (2.2.2), we reach to the
following ordinary dif\/ferential equation: 
\be \label{malek:I.54}
F^n F_{\eta\eta} +n F^{n-1}(F_\eta)^2-\sigma^2 F^s =-C_2 {\mathrm
e}^{-\xi\eta}.
\ee
For the case when $m=s=1$, and from (\ref{malek:I.51}), we f\/ind $n=0$.
Hence (\ref{malek:I.54}) becomes 
\be \label{malek:I.55}
F_{\eta\eta} -\sigma^2 F=-C_2 {\mathrm e}^{-\xi \eta},
\ee
and the corresponding boundary conditions are:
\be \label{malek:I.56}
\mbox{(i)} \ \quad   F(0) = 0,
\ee
\be \label{malek:I.57}
\mbox{(ii)}  \quad F(h)  =\gamma.
\ee
Applying boundary conditions (\ref{malek:I.56}) and (\ref{malek:I.57}), we get
the solution 
\be \label{malek:I.58}
F(\eta) =a_1{\mathrm e}^{\sigma \eta}+ a_2 {\mathrm e}^{-\sigma \eta}
+\frac{C_2}{\sigma^2-\xi^2} {\mathrm e}^{-\xi \eta}.
\ee
For f\/inite temperature, $a_1 = 0$, and applying boundary condition
(\ref{malek:I.56}), we get the solution in the form: 
\be \label{malek:I.59}
F(\eta) =\frac{C_2}{\sigma^2-\xi^2}
\left({\mathrm e}^{-\xi \eta} -\frac{\xi}{\sigma}{\mathrm e}^{-\sigma
\eta}\right).
\ee
Hence the temperature distribution across the lake, corresponding to
case (2) is: 
\be \label{malek:I.60}
T(z,t) =T_0+\frac{C_2 t}{\sigma^2-\xi^2} \left({\mathrm e}^{-\xi z}
-\frac{\xi}{\sigma} {\mathrm e}^{-\sigma z}\right).
\ee

For $0 < t < 150$ (in days), following Girgis and Smith [18], we use
the following values of the parameters: 
$C_2=0.2014$, $\beta = 12391$, $\xi = 0.048$, $h = 400$~meter, $T_0 =
4\,{}^\circ$C, the obtained results are plotted in Fig.3 and
Fig.4. 

\setcounter{equation}{0}

\section{$\!$Application (II): Unsteady free-convective boundary-layer
f\/low on a non-isothermal vertical f\/lat plate}

\subsection{Mathematical formulation}

Consider a natural convective, laminar, boundary layer adjacent to a
semi-inf\/inite, vertical f\/lat plate. The f\/luid is isothermal and of
constant temperature $\overline{T}_\infty$, far from the plate. The
plate has nonuniform surface temperature $\overline{T}_w>
\overline{T}_\infty$ (heated plate case). The f\/luid has the following
constant properties: ``$\beta$'' is the volumetric coef\/f\/icient of thermal
expansion, ``$\nu$'' the kinematic viscosity, and ``$\alpha$''is the
thermal dif\/fusivity. 

Along with the application of the Boussinesq and boundary-layer
approximation, the equations of motion may be written as:

(1)     Conservation of mass:
\be \label{malek:II.1}
\overline{u}_{\overline{x}}+ \overline{v}_{\overline{y}}=0.
\ee

(2)     Momentum equation:
\be \label{malek:II.2}
\overline{u}_{\overline{t}}+\overline{u}\,\overline{u}_{\overline{x}}+
\overline{v}\,\overline{u}_{\overline{y}}=g \beta(\overline{T}-
\overline{T}_\infty)+\nu \overline{u}_{\overline{y}\, \overline{y}}.
\ee

(3)     Energy equation: 
\be \label{malek:II.3}
\overline{T}_{\overline{t}}+\overline{u} \,\overline{T}_{\overline{x}}+
\overline{v}\, \overline{T}_{\overline{y}}=\alpha
\overline{T}_{\overline{y}\,\overline{y}}.
\ee

\setcounter{equation}{0}
\renewcommand{\theequation}{\arabic{section}.4.\arabic{equation}}

Boundary conditions:
\be \label{malek:II.4.1}
\mbox{(i)} \ \quad \overline{v}=0, \ \overline{u}=0, \
\overline{T}=\overline{T}_w(\overline{x}, \overline{t}\,) \quad
\mbox{at} \ \ \overline{y}=0, \ \overline{t}>0,
\ee
\be \label{malek:II.4.2}
\mbox{(ii)} \quad \overline{u}=0, \ \overline{T}=\overline{T}_\infty
\quad \mbox{as} \ \  \overline{y}\to \infty.
\ee

\setcounter{equation}{4}
\renewcommand{\theequation}{\arabic{section}.\arabic{equation}}

Dimensionalize the variables according to:
\be \label{malek:II.5}
\ba{l}
\ds x=\frac{\overline{x}}{L}, \qquad
y=(\mbox{Gr})^{1/4}\frac{\overline{y}}{L}, \qquad
T=\frac{\overline{T}-\overline{T}_\infty}{\Delta T},
\qquad
\Theta =\frac{T}{T_w}, \\[4mm]
\ds  u=\frac{\overline{u}}{U},
\qquad
v= (\mbox{Gr})^{1/4}\frac{\overline{v}}{U}, \qquad
t=U\frac{\overline{t}}{L},
\ea
\ee
where $L$ is some arbitrary reference length, $\Delta
T=\overline{T}_{\mbox{\scriptsize \rm ref}}-\overline{T}_\infty$,
$\overline{T}_{\mbox{\scriptsize \rm ref}}$    is
some arbitrary reference temperature, $U=(g\beta L \Delta T)^{1/2}$
is the Characteristic velocity, and 
$\ds \mbox{Gr}= \frac{g \beta L^3 \Delta T}{\nu^2}$ is the Grashof number.

\smallskip

\noindent
{\bf In dimensionalized form}

The basic equations are:
\be \label{malek:II.6}
u_x+v_y=0,
\ee
\be \label{malek:II.7}
u_t+u u_x +v u_y=T+u_{yy},
\ee
\be \label{malek:II.8}
T_t+uT_x +vT_y =\frac{1}{\mbox{Pr}} T_{yy}; \qquad Pr=\frac{\nu}{\alpha} \ 
\ \mbox{is the Prandtl number}.
\ee

\setcounter{equation}{0}
\renewcommand{\theequation}{\arabic{section}.9.\arabic{equation}}

The boundary conditions are:
\be \label{malek:II.9.1}
\mbox{(i)} \ \quad v=0, \ \ u=0, \ \ T=T_w(x,t) \ \ \mbox{at} \ \ y=0, \ t>0,
\ee
\be \label{malek:II.9.2}
\mbox{(ii)} \quad u=0, \ \ T=0, \ \  \mbox{as} \ \ y\to \infty.
\ee

\setcounter{equation}{9}
\renewcommand{\theequation}{\arabic{section}.\arabic{equation}}

From the continuity equation there exists a stream function $\Psi(x,y)$
such that 
\be \label{malek:II.10}
u=\frac{\p \Psi}{\p y}, \qquad v=-\frac{\p \Psi}{ \p x}.
\ee
Momentum and energy equations take the form:
\be \label{malek:II.11}
\Psi_{yt}+\Psi_y \Psi_{yx}-\Psi_x\Psi_{yy}=\Theta T_w +\Psi_{yyy},
\ee
\be \label{malek:II.12}
T_w \Theta_t +\Theta (T_w)_t +T_w \Psi_y \Theta_x +\Theta
\Psi_y(T_w)_x -T_w\Psi_x \Theta_y=\frac{1}{\mbox{Pr}} T_w\Theta_{yy}.
\ee
Boundary conditions are:
\be \label{malek:II.13}
\ba{l}
\ds \Psi_x(x,0,t)=\Psi_y(x,0,t)=0, \qquad \theta(x,0,t)=1,\\[3mm]
\ds \lim_{y\to \infty}\Psi_y(x,y,t)=0, 
\qquad \lim_{y\to \infty} \Theta(x,y,t)=0.
\ea
\ee

\subsection{Solution of the problem}

A class of two-parameter group $(a_1, a_2)$ has the form
\be \label{malek:II.14}
\overline{S}=C^S(a_1,a_2)S+K^s(a_1,a_2),
\ee
where ``$S$'' stands for $x$, $y$, $t$; $\psi$, $\Theta$, $T_w$, and
$C$'s and $K$'s are real valued and dif\/ferentiable functions with
respect to $a_1$ and $a_2$.

\smallskip

\noindent
{\bfseries \itshape The invariance analysis:}

Invariance of the transformed momentum equation leads to the
following group $G$: 
\be \label{malek:II.15}
G: \ \left\{ \ba{l}
\ds \overline{x} =(C^y C^\Psi) x+K^x, \quad \overline{y} =C^y y,
\quad \overline{t}=(C^y)^2t+K^t, \\[2mm]
\ds \overline{\Psi}=C^\Psi \Psi +K^\Psi, \quad
\overline{T}_w=\left(\frac{C^\psi}{(C^y)^3}\right) T_w, \quad
\overline{\Theta}=\Theta. 
\ea\right.
\ee

\noindent
{\bfseries \itshape The absolute invariants:}

A function $\eta(x,y,t)$ is an absolute invariant of a two-parameter
group if $\eta$ satisf\/ies: 
\be \label{malek:II.16}
\ba{l}
\ds (\alpha_1 x +\alpha_2) \eta_x +(\alpha_3 y +\alpha_4) \eta_y
+(\alpha_5 t+\alpha_6)\eta_t=0, \qquad \mbox{and}\\[2mm]
\ds (\beta_1 x+\beta_2 )\eta_x+(\beta_3 y+\beta_4) \eta_y+(\beta_5
t+\beta_6) \eta_t=0,
\ea
\ee
where
\[
\ba{l}
\ds \alpha_1 =\frac{\p C^x}{\p a_1}(a_1^0, a_2^0), \quad
\alpha_2= \frac{\p K^x}{\p a_1}(a_1^0, a_2^0), \
\ldots, \\[5mm]
\ds  \beta_1=\frac{\p C^x}{\p a_2}(a_1^0, a_2^0), \quad
\beta_2=\frac{\p K^x}{\p a_2}(a_1^0, a_2^0), \
\ldots .
\ea
\]

The only possible form for the absolute invariant $\eta$ is:
\be \label{malek:II.17}
\eta =\frac{y}{\sqrt{a_1 t +b_1}},
\ee
where $a_1 = \alpha_5 =\beta_5 $, $b_1 = \alpha_6 = \beta_6$ are
constants. 

Abd-el-Malek, {\it et al} [8] have shown that the case of 
$\eta=y/\sqrt{a_1x +b_1 t+c}$ does not lead to any solution.

\smallskip

\noindent
{\bfseries \itshape
The complete set of absolute invariants:}

The only possible forms for the absolute invariant $\Psi$ and $T_w$
are: 
\be \label{malek:II.18}
\Psi(x,y,t)=\Gamma(x,t) F(\eta),
\ee
\be \label{malek:II.19}
T_w=T_0 \,\omega(x,t).
\ee

\subsubsection{Solution corresponding to the form of $\eta$ in
(\ref{malek:II.17}):} 

        The corresponding dif\/ferential equations are:
\be \label{malek:II.20}
\ba{l}
\ds F_{\eta\eta\eta} +\left(\frac{a_1}{2} \eta +F\right) F_{\eta\eta}
-(F_\eta)^2 +a_1 F_\eta +\Theta=0,\\[4mm]
\ds \frac{1}{\mbox{Pr}} \Theta_{\eta\eta}+ \left(\frac{a_1}{2} \eta +F\right)
\Theta_\eta +(2a_1- F_\eta)\Theta=0
\ea
\ee
with the boundary conditions:
\be \label{malek:II.21}
F(0)=F_\eta(0) =0, \quad \Theta(0)=1, \quad F_\eta(\infty)=0, \quad
\Theta(\infty)=0. 
\ee
We get the following solution:
\be \label{malek:II.22}
\ba{l}
\ds T_w=\frac{x+b_2}{(0.4472 a_1 t+b_1)^2},
\qquad u=\frac{x+b_2}{a_1 t+b_1} F_\eta, \\[5mm]
\ds v =-\frac{F}{\sqrt{a_1 t+b_1}}, \qquad
q=(x+b_2) \frac{[-\Theta_\eta(0)]}{(a_1 t+b_2)^{5/2}}.
\ea
\ee

\subsection{Conclusion}

It is clear from the obtained results (\ref{malek:II.22}) that the temperature
prof\/ile overshoots in the region of the boundary layer near the
plate. This phenomena occurs for values of $a_1 > 1$, and becomes
stronger as $a_1$ increases. This means that this phenomenon is
accomplished by those cases for which $T_w$ decreases rapidly with
time. 

If we study the ef\/fect of Pr on the temperature prof\/ile we f\/ind
that there is a rapid increase in $\theta$ near the plate. 
This becomes more evident for larger values of Pr. Also the thermal
boundary-layer thickness decreases for increasing values of Pr.

\setcounter{equation}{0}

\section{Application (III): Dispersion of gaseous pollutants in the
presence of a temperature inversion}

\subsection{Mathematical formulation}

The gaseous pollutant is bounded from above by the ground surface and
from below by the inversion layer, 
which is at height ``$h$'' from the ground surface. Assuming that the
pollution, with concentration $C(x,y)$, is evenly 
distributed throughout the layer, and the mean concentration of the
pollutant at $x=0$ averaged over $0 \leq  y \leq  h$ is 
constant and equal to $C_0$. The dif\/fusion of the pollutants takes
place due to the wind that has a constant mean 
velocity $u = u(x)$ in the $x$-direction, and the eddy dif\/fusivities
$\kappa_1$ and $\kappa_2$ in the $x$ and $y$-directions,
respectively, are also independent of $y$.

The normalized steady state dif\/fusion equation, that governs the
dispersion of the gaseous pollutants is 
\be \label{malek:III.1}
u(x)C_x=\kappa_1(x) C_{xx}+\kappa_2(x) C_{yy},
\ee
with the boundary conditions
\be \label{malek:III.2}
\left.
\ba{ll}
\kappa_1C_y=\lambda \gamma C & \mbox{at} \ y=0\\[1mm]
C_y=0 & \mbox{at} \ y=1\\[1mm]
C=1 & \mbox{at} \ x=0, \ 0<y<1\\[1mm]
C_x=0 & \mbox{at} \ x\to \infty
\ea \right\}, \qquad \gamma^2=\frac{h}{u_0},
\ee
where all $x$ and $y$ are scaled with respect to $h$, $C$ with
respect to $C_0$, $u$ with respect to $u_0$, $\kappa_1$ and
$\kappa_2$ with respect to $u_0h$, and $u_0$ is a reference velocity.
Values of $\lambda$ classify two cases, case (1): $\lambda \ll 1$
corresponds to the case where no pollutant is absorbed by the ground,
case (2): $\lambda \gg 1$ corresponds to the case where all pollutant is
absorbed by the ground.

Introduce the non-dimensional function $\theta(x,y)$ and $C^*(x)$ such that
\be \label{malek:III.3}
C(x,y) = \theta(x,y) C^*(x), 
\ee
equation (\ref{malek:III.1}) becomes
\be \label{malek:III.4}
u(C^* \theta_x +\theta C_x^*)=\kappa_1(C^* \theta_{xx}+2\theta_x
C_x^*+\theta C^*_{xx})+\kappa_2C^*\theta_{yy}.
\ee

\subsection{Solution of the problem}

Following the same analysis as we did in application (I), we f\/ind
that the group $G_1$ which transforms invariantly 
the dif\/ferential equation (\ref{malek:III.4}) and the boundary conditions
(\ref{malek:III.2}), is in the form: 
\be \label{malek:III.5}
G_1: \quad \left\{
\ba{l}
\ds \overline{x}=E^x(a)x, \ \ \overline{y}=y, \ \
\overline{u}=E^x(a)u,\\[2mm] 
\ds \overline{\kappa}_1=(E^x(a))^2\kappa_1, \ \
\overline{\kappa}_2=\kappa_2,\ \
\overline{C}{\,}^*=E^{C^*}(a)C^*, \ \ \overline{\theta}=\theta.
\ea
\right.
\ee
The absolute invariants of the independent and dependent variables are:
\be \label{malek:III.6}
\ba{l}
\ds \eta=y, \qquad F(x,u)=\frac ux, \qquad
G(x,\kappa_1)=\frac{\kappa_1}{x^2}, \\[3mm]
\ds \theta=\theta(y), \qquad
\kappa_2=\kappa_2(x), \qquad C^*=C^*(x).
\ea
\ee

\subsection{The reduction to ordinary dif\/ferential equation}

Substituting from (\ref{malek:III.6}) into (\ref{malek:III.4}) gives:
\be \label{malek:III.7}
[\kappa_2 C^*]\theta_{yy} +[x^2GC_{xx}^*-x F C_x^*]\theta =0.
\ee
The requirement of reducing (\ref{malek:III.7}) to two ordinary dif\/ferential
equations, for some constant $p^2$, implies that: 
\be \label{malek:III.8}
[x^2GC_{xx}^* -x F C_x^*] =p^2[\kappa_2 C^*].
\ee
Under this assumption, (\ref{malek:III.4}) gives the ordinary dif\/ferential
equation of $\theta(y)$, namely
$\theta_{yy}+p^2 \theta =0$, which has the solution
\be \label{malek:III.9}
\theta(y) = A \cos p(y+\varepsilon).
\ee
Rearranging (\ref{malek:III.8}), we get
\be \label{malek:III.10}
\alpha C_{xx}^* -2\beta C_x^* -p^2 C^*=0,
\ee
where, using (\ref{malek:III.6}), the constants $\alpha$ and $\beta$ are
\be \label{malek:III.11}
\alpha=\frac{x^2 G}{\kappa_2}=\frac{\kappa_1}{\kappa_2} \quad
\mbox{and} \quad \beta=\frac{xF}{2\kappa_2}=\frac{u}{2\kappa_2}.
\ee
The ordinary dif\/ferential equation (\ref{malek:III.10}) has the general
solution 
\be \label{malek:III.12}
C^*(x)=A_1 {\mathrm e}^{m_1 x} +A_2 {\mathrm e}^{m_2 x}.
\ee
Applying the boundary conditions (\ref{malek:III.2}) as $x \to \infty$, for
$\alpha > 0$, we get  the solution 
\be \label{malek:III.13}
C^*(x)=A_1 {\mathrm e}^{m x},
\ee
where
\be \label{malek:III.14}
m=\beta-\sqrt{\beta^2+\alpha p^2}.
\ee
Substituting (\ref{malek:III.9}) and (\ref{malek:III.13}) in (\ref{malek:III.3}), we
get 
\be \label{malek:III.15}
C(x,y)=B {\mathrm e}^{m x}\cos p(y+\varepsilon).
\ee
Application of the boundary condition (\ref{malek:III.2}) at the inversion level
determines $\varepsilon = -1$. Hence 
\be \label{malek:III.16}
C(x,y)=B {\mathrm e}^{m x}\cos p(y-1).
\ee
Application of the boundary condition (\ref{malek:III.2}) at the ground surface
yields 
\be \label{malek:III.17}
\tan p=\frac{\lambda \gamma}{\kappa_2 p}.
\ee
The constants ``$B$'' and ``$p$'' will be determined corresponding to
the two limiting cases of~$\lambda$.

\subsection{Analytical solution corresponding to the two limiting
cases of $\lambda$}

{\bf Case (1):} $\lambda \ll 1$, which corresponds to no pollutant
absorbed by the ground. 

Assuming $\ds \frac{\lambda \gamma}{\kappa_2 p}$ to be a very small
quantity leads to 
\be \label{malek:III.18}
p^2=\frac{\lambda\gamma}{\kappa_2}
\ee
and hence, by applying the boundary condition (\ref{malek:III.2}), we f\/ind $B=1$.

The concentration distribution for the case (1) is
\be \label{malek:III.19}
C= {\mathrm e}^{m x}.
\ee
where ``$m$'' is given by (\ref{malek:III.14}).

\smallskip

\noindent
{\bf Case (2):} $\lambda \gg 1$, which corresponds to all pollutant 
absorbed by the ground. 

In this case, from equation (\ref{malek:III.17}), we get
\[
p=N \pi.
\]
The concentration distribution model will be
\be \label{malek:III.20}
C(x,y)=\sum_{n=1}^\infty B_n {\mathrm e}^{m x}\cos [N\pi(y-1)],
\ee
where
\be \label{malek:III.21}
m=\beta-\sqrt{\beta^2 +N^2 \pi^2 \alpha}, \qquad
N=\frac{2n-1}{2}; \quad n=1,2,3,\ldots\ .
\ee
The constants $B_n$ will be determined as the cosine's Fourier
coef\/f\/icients of the expansion of the function: $C=1$ in 
$ 0 < y < 1$ and $x = 0$, we get
\[
B_n=2\left( \frac{\sin N\pi}{N\pi}\right).
\]
Hence the concentration distribution for the case (2) is
\be \label{malek:III.22}
C(x,y)=2\sum_{n=1}^\infty \left( \frac{\sin N\pi}{N\pi}\right)
 {\mathrm e}^{m x}\cos [N\pi(y-1)],
\ee
where ``$m$'' and ``$N$'' are given by (\ref{malek:III.21}).

\subsection{Results and discussion}

For the case where no pollutant is absorbed by the ground surface, it
is found that the concentration distribution 
has the form
\[
C=\exp\left\{ \frac{1}{2\kappa_2}\left[ u-\sqrt{u^2+4\lambda
\kappa_1\gamma} \right]x\right\},
\]
which is independent of ``$y$'' and mainly depends on $u$,
$\kappa_1$, $\kappa_2$, $h$ and $\lambda$. For very small $\lambda$ or
very large $u$, it is clear that no pollutant will be absorbed by the
ground surface. Also we concluded that: 

\begin{enumerate}

\item[(1)] As ``$u$'' increases, the absorbed pollutant by the ground
surface will be less. 

\item[(2)] As the eddy dif\/fusivity ratio, $\alpha$, increases the absorbed
pollutant by the ground surface will be more. 
\end{enumerate}

For the case where all pollutant is absorbed by the ground surface,
it is found that the concentration distribution  has the form:
\[
C=2\sum_{n=1}^\infty \left( \frac{\sin N\pi}{N\pi}\right)
\exp\left[\frac{1}{2\kappa_2}\left( u-\sqrt{u^2+4N^2\pi^2 \kappa_1} 
\right)x\right] \cos [N\pi(y-1)].
\]

\end{document}